\documentclass[11pt,a4paper]{amsart}
\usepackage{amssymb}
\usepackage[pdfborder={0 0 0},pdfpagelayout=OneColumn]{hyperref}

\newtheorem{theorem}{Theorem}[section]
\newtheorem{lemma}[theorem]{Lemma}
\newtheorem{corollary}[theorem]{Corollary}
\newtheorem{proposition}[theorem]{Proposition}
\theoremstyle{remark}
\newtheorem{remark}[theorem]{Remark}
\theoremstyle{plain}
\newtheorem{geometricprinciple}[theorem]{Theorem}

\numberwithin{equation}{section}
\DeclareMathOperator{\tr}{tr}
\DeclareMathOperator{\codim}{codim}
\DeclareMathOperator{\dist}{dist}
\DeclareMathOperator{\Graph}{Graph}
\DeclareMathOperator{\spn}{span}
\newcommand{\Nuc}{\mathcal N}
\newcommand{\norm}[1]{\left\lVert#1\right\rVert}
\newcommand{\cl}[1]{\overline{#1}}
\newcommand{\R}{\mathbb R}
\newcommand{\C}{\mathbb C}

\begin{document}

\title[Hereditary Approximation under Hilbertian Extension]
{Hereditary Approximation under Hilbertian Extension}

\author{Gang Cai}
\address[G. Cai]{School of Mathematical Sciences, Chongqing Normal University, Chongqing 401331, China}
\email{caigang-aaaa@163.com}

\author{Qingjin Cheng}
\address[Q. Cheng]{School of Mathematical Sciences, Xiamen University, Xiamen 361005, China}
\email{qjcheng@xmu.edu.cn}

\author{Longfa Sun$^\dag$}
\address[L. Sun]{Hebei Key Laboratory of Physics and Energy Technology, School of Mathematics and Physics, North China Electric Power University, Baoding 071003, China}
\email{sun.longfa@ncepu.edu.cn}

\author{Yinghua Sun}
\address[Y. Sun]{School of Mathematical Sciences, Xiamen University, Xiamen 361005, China}
\email{sunak87@163.com}

\thanks{$^\dag$ Corresponding author.\\
This work is supported by the National Natural Science Foundation of China
(Grant Nos. 12571141 and 12071389), the Chongqing Talent Support Program
(Grant No. cstc2024ycjh-bgzxm0121), and the Fundamental Research Funds
for the Central Universities (Grant No. 2025MS178).
}

\date{}
\begin{abstract}
In this paper, we prove that if a Banach space $X$ has the hereditary approximation property, then
$X\oplus_2 H$ has the hereditary approximation property for every
Hilbert space $H$. This answers a question of Johnson and Szankowski.
In addition, if $X$ is also reflexive, then every closed subspace of $X\oplus_2H$
has the metric approximation property.
The proof uses Grothendieck's characterization of the approximation
property, a factorization of nuclear operators through quotients of $X$,
and a nuclear perturbation
of the coordinate projection on a closed graph.
\end{abstract}
\renewcommand{\subjclassname}{Mathematics Subject Classification (2020)}
\maketitle
\section{Introduction}
Let $X$ be a Banach space, it is said to have the \emph{approximation property}
(AP) if, for every compact subset $K$ of $X$ and every $\varepsilon>0$,
there is a bounded finite-rank operator $T:X\to X$ such that
$\norm{Tx-x}<\varepsilon$ for all $x\in K$. It has the
\emph{hereditary approximation property} (HAP) if every closed subspace
of $X$ has the AP. Hilbert spaces have the HAP, since their closed
subspaces are again Hilbert spaces. The corresponding assertion for
spaces with the AP is false: Szankowski constructed subspaces of
$\ell_p$, $1\le p<2$, which fail the AP~\cite{Szankowski}.

We shall also consider the bounded and metric approximation
properties. A Banach space $X$ has the \emph{bounded approximation
property} (BAP) if there is $\lambda\ge1$ such that the operators
$T$ in the definition of the AP can always be chosen with
$\norm{T}\le\lambda$. If they can be chosen with $\norm{T}\le1$,
then $X$ has the \emph{metric approximation property} (MAP).
We say that $X$ has the \emph{hereditary bounded approximation
property} (HBAP) if every closed subspace has the BAP, with a constant
that may depend on the subspace. Similarly, $X$ has the
\emph{hereditary metric approximation property} (HMAP) if every
closed subspace, with its inherited norm, has the MAP. Thus
\[
 \mathrm{HMAP}\Longrightarrow\mathrm{HBAP}\Longrightarrow\mathrm{HAP}.
\]
Hilbert spaces have the HMAP, since finite-dimensional orthogonal
projections give the MAP on each closed subspace.

The first examples of Banach spaces not isomorphic to Hilbert space
which have the HAP were constructed by Johnson~\cite{Johnson}.
These spaces are asymptotically Hilbertian. Recall that a space $X$
is \emph{asymptotically Hilbertian} if there is $\beta<\infty$ such
that, for every positive integer $n$, some finite-codimensional
subspace of $X$ has all its $n$-dimensional subspaces
$\beta$-isomorphic to $\ell_2^n$. Pisier subsequently proved that
weak Hilbert spaces have the HAP~\cite{Pisier}. As explained
in~\cite{JS}, both Johnson's examples and weak Hilbert spaces are
asymptotically Hilbertian, and an asymptotically Hilbertian space with
a symmetric basis must be isomorphic to $\ell_2$.

Johnson and Szankowski~\cite{JS} obtained a different sufficient
condition for the HAP. For an infinite-dimensional Banach space $X$, put
\[
 d_n(X)=\sup\{d(E,\ell_2^n):E\subset X,\ \dim E=n\},
\]
where $d$ denotes the Banach--Mazur distance. Their Theorem~2.1 gives a sufficient condition for the AP in terms of
sufficiently slow growth of $d_n(X)$. The condition passes to
subspaces and therefore implies the HAP. They also
constructed spaces with a symmetric basis, not isomorphic to $\ell_2$,
which satisfy this hypothesis. Thus the HAP does not imply that a
space is asymptotically Hilbertian. Further examples were obtained by
Morales and Su\'arez de la Fuente~\cite{MS}, who showed that the
twisted Hilbert space $Z(T_s^2)$ has the HAP but is not asymptotically
Hilbertian.

For reflexive spaces, the HAP, the HBAP, and the HMAP are equivalent.
Indeed, a closed subspace of a reflexive HAP space is reflexive and has the
AP, and Grothendieck's theorem gives the MAP for every reflexive
space with the AP (see~\cite[Theorem~1.e.15 and pp.~40--41]{LT}).
The implication from the AP to the BAP in the reflexive case is also
used in~\cite[p.~1998]{JS}. Thus Johnson's HAP
examples and weak Hilbert spaces have the HMAP, since they are
superreflexive. The spaces with a symmetric basis constructed
in~\cite{JS} also have the HMAP. Their superreflexivity follows from
the slow growth of their local Euclidean distances, as recalled
in~\cite[p.~1992]{JS}. These give examples of HMAP spaces not
isomorphic to Hilbert space.

Although these results give several classes of spaces with the HAP,
the property itself has strong geometric consequences. As recalled
in~\cite[p.~1988]{JS}, an HAP space has type $p$ for every $1<p<2$
and cotype $q$ for every $2<q<\infty$. A space with both type~2 and
cotype~2 is isomorphic to a Hilbert space, by Kwapie\'n's
theorem~\cite{Kwapien}. The estimates just stated do not give these
endpoint properties. In particular, they do not reduce questions
about the HAP to questions about Hilbert spaces.

One difficulty in working with the HAP is its behavior under direct
sums. The AP passes to finite direct sums, whereas the HAP need not:
Casazza, Garc\'ia, and Johnson~\cite{CGJ} constructed two spaces with
the HAP whose direct sum fails the HAP. See also~\cite[p.~1989]{JS}.
This leaves the question of what happens when one summand is Hilbert.
Johnson and Szankowski asked the following question
in~\cite[Question~5, p.~1999]{JS}: if $X$ has the HAP, must
$X\oplus\ell_2$ have the HAP? The choice of an equivalent norm on
this finite direct sum does not affect the AP or the HAP. We use
the $\ell_2$ direct sum norm and prove the following theorem, which
gives an affirmative answer to their question.

\begin{theorem}\label{thm:main}
If a Banach space $X$ has the hereditary approximation property, then
$X\oplus_2\ell_2$ has the hereditary approximation property.
\end{theorem}

By a separable reduction, the same conclusion holds for an arbitrary
Hilbert summand, as stated in Corollary~\ref{cor:arbitrary-hilbert}.
If $X$ is reflexive, the resulting direct sum is reflexive as well.
The equivalence just noted therefore gives the HMAP for
$X\oplus_2H$. We record this consequence in
Corollary~\ref{cor:hmap}. In particular, the HMAP, and hence the
HBAP, is preserved under the addition of a Hilbert summand within
the class of reflexive Banach spaces.

However, the above theorem does not resolve the question about $\ell_2(X)$
in~\cite[Question~9]{JS}. It adds a Hilbert space to $X$, whereas
$\ell_2(X)$ involves countably many copies of $X$. Related results
on subspaces of $\ell_2(X)$ without the AP were obtained by Anisca
and Chlebovec~\cite{AC}. We return to this distinction after the
proof of the theorem.

Let us describe the proof for a separable Hilbert space $H$. The
case of an arbitrary Hilbert space follows by separable reduction.
If $Y$ is a closed subspace of $X\oplus_2H$, its projection onto
$X$ need not be closed.
Consequently, even though $X\oplus_2H$ has the AP, the AP of $Y$
cannot be obtained simply by decomposing $Y$ along the two summands.
We use instead Grothendieck's characterization of the AP in terms of
the nuclear trace. If $Y$ fails the AP, it gives a nuclear operator
$T:X\oplus_2H\to Y$ which vanishes on $Y$ and satisfies
$\tr_{X\oplus_2H}(\iota T)\ne0$, where $\iota:Y\to X\oplus_2H$
is the inclusion. Writing
\[
 T(x,h)=Sx+Rh,
\]
we obtain nuclear maps $S:X\to Y$ and $R:H\to Y$. If $J:Y\to X$
and $Q:Y\to H$ are the coordinate maps, their relation is
$SJ+RQ=0$.

There are two ingredients in the proof. The first applies when the
projection of $Y$ onto $X$ is contained in a closed subspace
$M\subset X$ of infinite codimension. With $q:X\to X/M$ the quotient
map, we construct a compact map $a:H\to X$ and nuclear maps
$V_0,V_1:X/M\to Y$ such that
\[
 V_0qa=R,\qquad V_1qa=Sa.
\]
The finite-dimensional Euclidean subspaces needed for the construction
are supplied by Theorem~\ref{thm:geometric} and local reflexivity.
The former combines the theorem of Figiel and Tomczak-Jaegermann
with Pisier's equivalence of B-convexity and $K$-convexity. The additional
point is to choose the same map $a$ for both factorizations, while
keeping $V_0$ and $V_1$ nuclear into $Y$.

To see why the second identity is needed, put $u=J+aQ$. The first
identity alone gives
\[
 (S+V_0q)u=SaQ,
\]
because $qJ=0$ and $SJ+RQ=0$. The second identity removes the term
on the right and gives
\[
 \bigl(S+(V_0-V_1)q\bigr)u=0.
\]
We can then construct a nuclear operator into $\cl{u(Y)}\subset X$
which vanishes on that subspace and whose composition with the
inclusion into $X$ has trace $\tr_{X\oplus_2H}(\iota T)$. This
contradicts the HAP of $X$. No AP assumption on $Y$ or $X/M$ is
made in either factorization.

The second ingredient reduces the general case to this one. After
passing to a separable subspace, replacing $X$ by the closure of
the horizontal projection, and splitting off the vertical Hilbert
subspace, a possible counterexample is the graph of a closed,
densely defined, unbounded operator $G:D(G)\subset X\to H$. We
construct a nuclear map $C:H\to X$, with arbitrarily small nuclear
norm, such that
\[
 \{x-CGx:x\in D(G)\}
\]
is contained in a closed subspace of $X$ of infinite codimension.
The isomorphism $(x,h)\mapsto(x-Ch,h)$ then puts the graph in the
position required for the first part of the proof. This construction
is proved for an arbitrary reflexive codomain. Both the isomorphism
and its inverse can be chosen arbitrarily close to the identity.
Together, the two constructions allow the HAP of $X$ to be applied
to an arbitrary closed subspace of $X\oplus_2H$.

This paper is organized as follows. In Section~\ref{sec:trace}, we recall Grothendieck's characterization and the separable reduction. In Sections~\ref{sec:sections} and~\ref{sec:factorization}, we present the finite-dimensional construction and the
two factorization identities. They are
applied to the trace in Section~\ref{sec:transfer}. The change of
coordinates for a closed graph is constructed in
Section~\ref{sec:shear}, and the proof of the theorem is completed in
Section~\ref{sec:completion}.

\section{Preliminaries}\label{sec:trace}

Throughout the paper, all Banach spaces are over a fixed scalar
field $\mathbb F$, either $\R$ or $\C$. All subspaces, linear spans,
and operators are understood over $\mathbb F$, and dimensions are
taken over this field. For a Banach space $E$, the dual space $E^*$
consists of continuous $\mathbb F$-linear functionals on $E$.
In the complex case, Hilbert space inner products are taken to be
linear in the first variable. Subspaces are understood to be closed
unless otherwise specified, and carry the norm of the ambient space. For Banach spaces
$E$ and $F$, we write $E\oplus_2F$ for the direct sum with norm
\[
 \norm{(e,f)}=(\norm{e}^2+\norm{f}^2)^{1/2}.
\]
The identity operator on $E$ is denoted by $I_E$. We omit the
subscript when the underlying space is clear.
The AP, the HAP, and the HBAP are invariant under isomorphisms.
In particular, the choice of an equivalent norm on a finite direct
sum does not affect any of these properties. For isomorphic Banach spaces $E,F$, their
\emph{Banach--Mazur distance} is
\[
 d(E,F)=\inf\{\norm{T}\norm{T^{-1}}:
 T:E\to F\text{ is an isomorphism}\}.
\]
We use this notation only for finite-dimensional spaces of the same
dimension. A finite-dimensional space $E$ is $K$-Euclidean if
$d(E,\ell_2^{\dim E})\le K$. A subspace $F\subset E$ is
$C$-complemented in $E$ if there is a bounded projection $P:E\to F$
with $P|_F=I_F$ and $\norm{P}\le C$.

For $1\le p\le2$, a Banach space $E$ has \emph{type $p$} if there
is $T_p<\infty$ such that every finite family $x_1,\ldots,x_n\in E$
satisfies
\[
 \left(\mathbb E\left\lVert\sum_{j=1}^n\varepsilon_jx_j
 \right\rVert^2\right)^{1/2}
 \le T_p\left(\sum_{j=1}^n\norm{x_j}^p\right)^{1/p},
\]
where the $\varepsilon_j$ are independent random signs taking the
values $1$ and $-1$ with equal probability. For $2\le q<\infty$, it
has \emph{cotype $q$} if there is $C_q<\infty$ such that
\[
 \left(\sum_{j=1}^n\norm{x_j}^q\right)^{1/q}
 \le C_q\left(\mathbb E\left\lVert\sum_{j=1}^n\varepsilon_jx_j
 \right\rVert^2\right)^{1/2}
\]
for every such family.

A Banach space $E$ is \emph{B-convex} if there are an integer $n\ge2$
and $\delta>0$ such that, for all $x_1,\ldots,x_n$ in its closed unit
ball, some choice of signs satisfies
\[
 \left\lVert\sum_{j=1}^n\varepsilon_jx_j\right\rVert
 \le (1-\delta)n.
\]
Equivalently, $E$ does not contain the spaces $\ell_1^n$ uniformly:
there is no $K<\infty$ such that, for every $n$, some subspace
$E_n\subset E$ satisfies $d(E_n,\ell_1^n)\le K$. B-convexity is also
equivalent to having type $p$ for some $p>1$, and to
$K$-convexity~\cite{PisierK}. We use the equivalence with
$K$-convexity only to apply the geometric theorem in
Section~\ref{sec:sections}.

For a subspace $N\subset E$, we write
$N^\perp=\{f\in E^*:f|_N=0\}$ for its annihilator and
$\codim_E N=\dim(E/N)$ for its codimension.

For Banach spaces $E,F$, let $\Nuc(E,F)$ denote the Banach space of
nuclear operators from $E$ to $F$, equipped with the nuclear norm $\nu$.
Thus $T\in\Nuc(E,F)$ has a representation
\[
 Tx=\sum_{k=1}^{\infty}f_k(x)y_k,
 \qquad f_k\in E^*,\quad y_k\in F,\quad
 \sum_{k=1}^{\infty}\norm{f_k}\norm{y_k}<\infty,
\]
and $\nu(T)$ is the infimum of these sums. We write $f\otimes y$
for the rank-one operator $x\mapsto f(x)y$. We shall always specify
the codomain of a nuclear operator. The
inclusion $T(E)\subset F_0\subset F$ alone does not imply that
$T:E\to F_0$ is nuclear. We use the ideal inequality
\begin{equation}\label{eq:ideal}
 \nu(BTA)\le\norm{B}\nu(T)\norm{A}.
\end{equation}
Completeness in the nuclear norm and~\eqref{eq:ideal} do not require
the AP of either space.

If $E$ has the AP, the trace of a nuclear operator on $E$ is well defined by
$\tr_E T=\sum_k f_k(y_k)$. Thus the sum is independent of the nuclear representation.
If $E$ and $F$ both have the AP,
$T\in\Nuc(E,F)$, and $B:F\to E$ is bounded, then
\begin{equation}\label{eq:cyclic}
 \tr_E(BT)=\tr_F(TB).
\end{equation}
Indeed, writing $T=\sum_k f_k\otimes y_k$, the maps $BT$ and $TB$ have
the nuclear representations $\sum_k f_k\otimes By_k$ and
$\sum_k(f_k\circ B)\otimes y_k$, respectively. Both traces therefore
equal the absolutely convergent scalar series $\sum_k f_k(By_k)$.

We use two facts recalled in~\cite[p.~1988]{JS}. First, a space with
the HAP has type $p$ for every $1<p<2$, and hence is B-convex. Second,
Grothendieck~\cite{Grothendieck} proved the following characterization
of the AP (see also~\cite[p.~1988]{JS}): a Banach space $Y$ fails
the AP if and only if there are $y_k^*\in Y^*$ and $y_k\in Y$
such that
\begin{equation}\label{eq:witness}
 \sum_k\norm{y_k^*}\norm{y_k}<\infty,\qquad
 \sum_k y_k^*(y)y_k=0\quad(y\in Y),\qquad
 \sum_k y_k^*(y_k)=1.
\end{equation}
The last sum can be normalized to~1 by rescaling the vectors.
Notice that these identities concern a representation of the zero
operator. No trace on a space without the AP is being used.

As in the discussion of the trace in~\cite[pp.~1988--1989]{JS},
Hahn--Banach extensions give the following useful consequence.

\begin{lemma}\label{lem:trace-witness}
Let $E$ have the AP, let $Y\subset E$ fail the AP, and let $\iota:Y\to E$ be
the inclusion. There is $T\in\Nuc(E,Y)$ such that
\[
 T\iota=0,\qquad \tr_E(\iota T)=1.
\]
\end{lemma}

\begin{proof}
Let $(y_k^*,y_k)$ satisfy~\eqref{eq:witness}, and extend $y_k^*$ to
$\widetilde y_k^*\in E^*$ without increasing its norm.
Then
\[
 Tx=\sum_k\widetilde y_k^*(x)y_k
\]
defines an element of $\Nuc(E,Y)$. By~\eqref{eq:witness},
$T\iota=0$. Since $E$ has the AP, its nuclear trace is well defined,
and
\[
 \tr_E(\iota T)=\sum_k\widetilde y_k^*(y_k)
 =\sum_k y_k^*(y_k)=1.
\]
\end{proof}

Grothendieck's result also shows that every space without the AP
contains a separable subspace without the AP.

\begin{lemma}\label{lem:separable-reduction}
If a Banach space $Y$ fails the AP, it contains a separable closed subspace
that fails the AP.
\end{lemma}

\begin{proof}
Choose $(y_k^*,y_k)$ as in~\eqref{eq:witness} and put
$Y_s=\cl{\spn\{y_k:k\ge1\}}$. The pairs
$(y_k^*|_{Y_s},y_k)$ satisfy all three identities
in~\eqref{eq:witness}. Hence $Y_s$ is a separable closed subspace
without the AP.
\end{proof}

\section{Finite-dimensional sections of quotient maps}\label{sec:sections}

The geometric input is the theorem of Figiel and
Tomczak-Jaegermann~\cite{FTJ} on complemented Euclidean subspaces,
together with Pisier's equivalence of B-convexity and
$K$-convexity~\cite{PisierK}. We recall the resulting statement in
the form recorded in Plichko's discussion of complemented Euclidean
subspaces~\cite{Plichko}.

\begin{geometricprinciple}[Figiel--Tomczak-Jaegermann and Pisier]\label{thm:geometric}
Let $Z$ be B-convex. There is $d_Z\ge1$ such that for every positive
integer $n$ and every $\eta>0$ there is a positive integer
$m=m(Z,n,\eta)\ge n$ with the following property. Every $m$-dimensional
subspace $E\subset Z$ contains an $n$-dimensional subspace $F$ such that
$d(F,\ell_2^n)\le1+\eta$ and there is a bounded linear map $R:Z\to F$ with
$R|_F=I_F$ and $\norm{R}\le d_Z$.
\end{geometricprinciple}

The fact that the projection acts on all of $Z$, with a bound
independent of the dimension, will be used below. Recall also that
the dual of a B-convex space is B-convex. The other fact needed is
the principle of local reflexivity, in the following form: if
$E_0\subset X^{**}$ and $F_0\subset X^*$ are finite dimensional and
$\delta>0$, there is a linear map $T:E_0\to X$ with
$\norm{T}\le1+\delta$ such that
\[
 f(T\xi)=\xi(f)\qquad(\xi\in E_0,\ f\in F_0).
\]
See~\cite[Lemma~5 and its proof]{Plichko} for its use in this setting.
Neither local reflexivity nor the preceding duality assertion requires
reflexivity of $X$.

The next lemma is a consequence of these standard facts. It is the
form of the local reflexivity argument that we shall use, with the
projection required to vanish on a given subspace of infinite
codimension. The passage from the dual to the original space is
closely related to~\cite[Lemma~5]{Plichko}. We include the argument
to specify the annihilation condition and the constants.

\begin{lemma}\label{lem:projection}
Let $X$ be B-convex. There is $D\ge1$, depending only on $X$, such that
for every subspace $N\subset X$ of infinite codimension and every
positive integer $n$ there is a rank-$n$ projection $P:X\to X$ satisfying
\begin{equation}\label{eq:projection}
 PN=0,\qquad\norm{P}\le D,\qquad d(PX,\ell_2^n)\le2D.
\end{equation}
\end{lemma}

\begin{proof}
Apply Theorem~\ref{thm:geometric} to $Z=X^*$ and put $d_0=d_{X^*}$.
Since $N^\perp$ is isometric to $(X/N)^*$, it is infinite dimensional. Choose in $N^\perp$ a subspace of dimension
$m(X^*,n,1)$. The theorem gives $F\subset N^\perp$, $\dim F=n$, and
$R_0:X^*\to F$ with
\[
 d(F,\ell_2^n)\le2,\qquad R_0|_F=I_F,\qquad\norm{R_0}\le d_0.
\]
In particular, $d_0$ is independent of $N$ and $n$.

Set $E_0=R_0^*(F^*)\subset X^{**}$. Local reflexivity with
$F_0=F$ and $\delta=1$ gives $T:E_0\to X$ such that $\norm{T}\le2$
and the pairings with $F$ are preserved. Let $s:X\to F^*$ be the
evaluation map, $s(x)(f)=f(x)$, and put $A=TR_0^*:F^*\to X$.
For $\varphi\in F^*$ and $f\in F$,
\[
 (sA\varphi)(f)=f(TR_0^*\varphi)
 =(R_0^*\varphi)(f)=\varphi(R_0f)=\varphi(f).
\]
Thus $sA=I_{F^*}$. In particular, $A$ is one-to-one and
$P=As$ satisfies
$P^2=AsAs=A(sA)s=As=P$. Its range is $A(F^*)$ and has dimension
$n$. Since $F\subset N^\perp$, $sN=0$, so $PN=0$.
Moreover,
\[
 \norm{P}\le2d_0,\qquad\norm{A}\le2d_0,
 \qquad \norm{A^{-1}:A(F^*)\to F^*}\le1,
\]
the last estimate following from $A^{-1}=s|_{A(F^*)}$ and
$\norm{s}\le1$. Finite-dimensional duality gives
$d(F^*,\ell_2^n)=d(F,\ell_2^n)$. Composing an almost optimal
isomorphism $\ell_2^n\to F^*$ with $A$ therefore yields
\[
 d(PX,\ell_2^n)\le2d_0\,d(F^*,\ell_2^n)\le4d_0.
\]
Taking $D=2d_0$ proves the assertion.
\end{proof}

For the factorization below, the projection lemma must be used with
two further conditions. We need the range of a section to lie in a
prescribed finite-codimensional subspace, and the corresponding map
on the quotient to vanish on a prescribed finite-dimensional
subspace. The next lemma imposes both conditions with uniform
bounds. In particular, the bounds do not depend on the codimension
of the constraint on the range.

\begin{lemma}\label{lem:sections}
Let $X$ be B-convex, let $M\subset X$ have infinite codimension, and
let $q:X\to X/M$ be the quotient map. There is a constant $C$,
depending only on $X$, with the following property. Given a
finite-dimensional subspace $F\subset X$, a finite-codimensional
subspace $L\subset X$, and a finite-dimensional Hilbert space $K$,
there are operators $t:K\to X$ and $\beta:X/M\to K$ such that
\begin{equation}\label{eq:sections}
 \begin{gathered}
 t(K)\subset L,\qquad\norm{t}\le1,\qquad\norm{\beta}\le C,\\
 \beta qt=I_K,\qquad (\beta q)|_F=0.
 \end{gathered}
\end{equation}
\end{lemma}

\begin{proof}
The case $K=\{0\}$ is immediate. Write $d=\dim K$ and
$c=\codim_X L$. The space $N=M+F$ is closed and has infinite codimension: indeed,
$qF$ is finite dimensional and $X/(M+F)$ is isomorphic to $(X/M)/qF$. Apply
Lemma~\ref{lem:projection} to $N$ with rank $d+c$. Since the
Banach--Mazur distance between finite-dimensional spaces is attained,
we may choose an isomorphism and rescale it to obtain
\[
 r:\ell_2^{d+c}\longrightarrow PX,
 \qquad\norm{r}\le1,\qquad\norm{r^{-1}}\le2D.
\]
Let $\pi_L:X\to X/L$ be the quotient map. Since
$\operatorname{rank}(\pi_Lr)\le c$, we have
\[
 \dim r^{-1}(L)=\dim\ker(\pi_Lr)\ge(d+c)-c=d.
\]
Choose an isometric embedding $j:K\to\ell_2^{d+c}$ with
$j(K)\subset r^{-1}(L)$, and define
\[
 t=rj,\qquad \beta q=j^*r^{-1}P,
\]
where $j^*$ is the Hilbert space adjoint. By the choice of $j$,
$t(K)\subset L$, and $\norm{t}\le\norm{r}\norm{j}\le1$. The composite
$j^*r^{-1}P:X\to K$ is bounded and vanishes on $M$ because $PM=0$. Hence it factors uniquely through $q$ to define $\beta:X/M\to K$.
The map $\beta$ annihilates $qF$ since $PF=0$,
and
\[
 \beta qt=j^*r^{-1}Prj=j^*j=I_K.
\]
Because $q$ is a norm-one quotient map, the quotient norm gives
$\norm{\beta}=\norm{\beta q}\le\norm{j^*}\norm{r^{-1}}\norm{P}
\le2D^2$.
Thus we may take $C=2D^2$.
\end{proof}

Applying the lemma successively gives mutually biorthogonal sections.
The finite-codimensional restrictions allow the new range to lie in
the kernels of all previously chosen coordinate maps.

\begin{corollary}\label{cor:blocks}
With $X,M,q,C$ as above, let $(K_i)$ be a sequence of finite-dimensional
Hilbert spaces and $(L_i)$ a sequence of finite-codimensional subspaces
of $X$. There are maps $t_i:K_i\to X$ and $\beta_i:X/M\to K_i$ such
that
\begin{equation}\label{eq:blocks}
 \norm{t_i}\le1,\quad\norm{\beta_i}\le C,\quad t_i(K_i)\subset L_i,
 \quad
 \beta_iqt_j=\begin{cases}I_{K_i},&i=j,\\0,&i\ne j.\end{cases}
\end{equation}
\end{corollary}

\begin{proof}
At stage $i$, apply Lemma~\ref{lem:sections} with
\[
 F=\sum_{j<i}t_j(K_j),\qquad
 L=L_i\cap\bigcap_{j<i}\ker(\beta_jq).
\]
Here $F$ is finite dimensional and $L$ has finite codimension.
The two choices give $\beta_iqt_j=0$ and $\beta_jqt_i=0$ for $j<i$,
while $\beta_iqt_i=I_{K_i}$ follows from~\eqref{eq:sections}.
The bounds are independent of $i$.
\end{proof}

\begin{remark}\label{rem:unweighted}
Neither convergence nor boundedness of the unweighted sum
$\sum_i t_i\beta_iq$ is asserted or needed. All infinite sums used
below have explicit summability estimates.
\end{remark}

\section{Simultaneous nuclear factorization}\label{sec:factorization}

Factorization of compact operators through auxiliary Banach spaces
is a classical tool in approximation theory. For example,
Johnson~\cite{JohnsonCompact} studied factorizations with compact
factors through a common separable space. Here the intermediate
space is the prescribed quotient $X/M$, and the maps from that
quotient must be nuclear into the prescribed space $Y$.

There is also a second requirement. As explained in the introduction,
a factorization of $R$ alone leaves an additional term in the
identity needed for the trace argument. We therefore construct, with
the same compact map $a:H\to X$, a nuclear factorization of $Sa$.
The finite-codimensional restrictions in Lemma~\ref{lem:sections}
allow the sections to be chosen in the kernels of finite-rank
approximations to $S$. This gives summability of their images under
$S$ in nuclear norm.

\begin{lemma}\label{lem:factorization}
Let $X$ be B-convex, let $M\subset X$ have infinite codimension, and
let $q:X\to X/M$ be the quotient map. Let $H$ be a separable Hilbert
space and $Y$ a Banach space. For every $S\in\Nuc(X,Y)$ and
$R\in\Nuc(H,Y)$, there are a compact operator $a:H\to X$ and
operators $V_0,V_1\in\Nuc(X/M,Y)$ such that
\begin{equation}\label{eq:factorization}
 V_0qa=R,\qquad V_1qa=Sa.
\end{equation}
\end{lemma}

\begin{proof}
If $H=\{0\}$, take all three maps to be zero. Otherwise choose
increasing finite-rank orthogonal projections $P_i$ on $H$ such that
$P_i\to I_H$ strongly and
\begin{equation}\label{eq:tail}
 \nu\bigl(R(I_H-P_i)\bigr)\le2^{-3i}\qquad(i\ge1).
\end{equation}
Indeed, fix a nuclear representation $R=\sum_k h_k^*\otimes y_k$ with
$\sum_k\norm{h_k^*}\norm{y_k}<\infty$. Start with orthogonal
projections onto increasing finite-dimensional spans of a fixed
orthonormal basis. For each $k$, the Riesz representation theorem
gives $u_k\in H$ such that $h_k^*(h)=\langle h,u_k\rangle$ for all
$h\in H$. Since $P$ is an orthogonal projection,
\[
 \norm{h_k^*(I_H-P)}=\norm{(I_H-P)u_k}\longrightarrow0
\]
as $P\to I_H$ strongly. Moreover,
$\norm{h_k^*(I_H-P)}\le\norm{h_k^*}$, so dominated convergence yields
\[
 \nu\bigl(R(I_H-P)\bigr)
 \le\sum_k\norm{h_k^*(I_H-P)}\norm{y_k}\longrightarrow0.
\]
A subsequence therefore satisfies~\eqref{eq:tail} and still
converges strongly to $I_H$. If $H$ is finite dimensional, take
$P_i=I_H$ for all sufficiently large $i$.

Set $P_0=0$, $E_i=P_i-P_{i-1}$, $H_i=E_iH$, and $\alpha_i=2^{-i}$.
We also regard $E_i$ as a map from $H$ onto $H_i$ and use zero
maps for zero blocks. For $i\ge2$, $H_i\subset(I_H-P_{i-1})H$.
Hence~\eqref{eq:tail} gives
$\nu(R|_{H_i})\le2^{-3(i-1)}$. It follows that
\begin{equation}\label{eq:weighted-tail}
 \sum_i\alpha_i^{-1}\nu(R|_{H_i})
 \le2\nu(R)+\sum_{i=2}^{\infty}2^{3-2i}
 =2\nu(R)+\frac23<\infty.
\end{equation}

Fix a nuclear representation $S=\sum_k x_k^*\otimes v_k$ with
$\sum_k\norm{x_k^*}\norm{v_k}<\infty$. For each $i$, truncate far enough to
obtain a finite-rank partial sum $S_i$ with
$\nu(S-S_i)\le2^{-i}$. The space
$L_i=\ker S_i$ is closed and finite codimensional. Apply
Corollary~\ref{cor:blocks} with $K_i=H_i$ and these $L_i$. The resulting
$t_i,\beta_i$ satisfy~\eqref{eq:blocks}, and
\begin{equation}\label{eq:S-small}
 \nu(St_i)=\nu((S-S_i)t_i)\le2^{-i}.
\end{equation}
Define
\begin{equation}\label{eq:factor-series}
 \begin{aligned}
 a&=\sum_i\alpha_it_iE_i,\\
 V_0&=\sum_i\alpha_i^{-1}(R|_{H_i})\beta_i,\\
 V_1&=\sum_i(St_i)\beta_i.
 \end{aligned}
\end{equation}
The series for $a$ converges absolutely in operator norm because
$\sum_i\norm{\alpha_it_iE_i}\le\sum_i\alpha_i=1$, and each term
has finite rank. Hence $a$ is a norm limit of finite-rank maps and is
compact. The other two series converge absolutely in $\Nuc(X/M,Y)$:
\[
 \begin{aligned}
 \sum_i\nu\bigl(\alpha_i^{-1}(R|_{H_i})\beta_i\bigr)
 &\le C\left(2\nu(R)+\frac23\right),\\
 \sum_i\nu\bigl((St_i)\beta_i\bigr)&\le C\sum_i2^{-i}=C.
 \end{aligned}
\]
In particular, $V_0,V_1\in\Nuc(X/M,Y)$.

By the operator-norm convergence of $a$ and biorthogonality, for each
$i$,
\[
 \beta_iqa=\sum_j\alpha_j\beta_iqt_jE_j=\alpha_iE_i.
\]
We may therefore compose the nuclear-norm convergent series for
$V_0$ and $V_1$ with $qa$ to obtain
\[
 V_0qa=\sum_i(R|_{H_i})E_i=R,
 \qquad
 V_1qa=\sum_i\alpha_iSt_iE_i=Sa.
\]
The partial sums in the first identity are $RP_m$, which converge
to $R$ in nuclear norm by~\eqref{eq:tail}. For the second identity,
the ideal inequality gives
\[
 \nu\left(Sa-\sum_{i=1}^m\alpha_iSt_iE_i\right)
 \le \nu(S)\left\|a-\sum_{i=1}^m\alpha_it_iE_i\right\|
 \longrightarrow0.
\]
Thus both identities follow from convergent series, with nuclearity
in $Y$ throughout the construction.
\end{proof}

\section{Nuclear operators and the hereditary approximation property}\label{sec:transfer}

We now apply the two factorization identities to nuclear operators
satisfying the relation $SJ+RQ=0$. The HAP of $X$ is used to
compute a trace on the closed subspace $\cl{(J+aQ)(Y)}$.
It is useful to state this step for arbitrary bounded maps $J$ and
$Q$. In the application they will be the coordinate maps of a
subspace of $X\oplus_2H$.

\begin{proposition}\label{prop:transfer}
Let $X$ have the HAP, let $H$ be a separable Hilbert space, and let
$Y$ be a Banach space. Let $J:Y\to X$ and $Q:Y\to H$ be bounded
linear maps. Suppose that
\begin{equation}\label{eq:horizontal}
 J(Y)\subset M
\end{equation}
for a subspace $M\subset X$ of infinite codimension. If
$S\in\Nuc(X,Y)$ and $R\in\Nuc(H,Y)$ satisfy
\begin{equation}\label{eq:compatible}
 SJ+RQ=0,
\end{equation}
then
\begin{equation}\label{eq:trace-zero}
 \tr_X(JS)+\tr_H(QR)=0.
\end{equation}
\end{proposition}

\begin{proof}
Since $X$ has the HAP, it is B-convex. Let $q:X\to X/M$ be the
quotient map. Apply Lemma~\ref{lem:factorization} to $q$, $S$, and
$R$, and set
\begin{equation}\label{eq:transfer-maps}
 V=V_0-V_1,\qquad u=J+aQ:Y\to X,\qquad A=S+Vq:X\to Y.
\end{equation}
Then $A\in\Nuc(X,Y)$, $qJ=0$, and $Vqa=R-Sa$. Consequently,
\begin{equation}\label{eq:Au}
 Au=SJ+SaQ+VqJ+VqaQ=SJ+RQ=0.
\end{equation}
Let $W=\cl{u(Y)}\subset X$, let $\iota_W:W\to X$ be the inclusion,
and write $u_W:Y\to W$ for $u$ with its codomain restricted to $W$.
Since $u_W$ is bounded, the operator
\[
 B=u_WA:X\to W
\]
is nuclear into $W$ by~\eqref{eq:ideal}. More explicitly, if
$A=\sum_k a_k^*\otimes y_k$ is a nuclear representation into $Y$,
then
\[
 B=\sum_k a_k^*\otimes u_Wy_k,\qquad
 \sum_k\norm{a_k^*}\norm{u_Wy_k}
 \le\norm{u_W}\sum_k\norm{a_k^*}\norm{y_k}<\infty.
\]
Thus nuclearity into $W$ follows from bounded composition, rather
than from restricting the codomain of an arbitrary nuclear map.
Equation~\eqref{eq:Au} implies that $B$ vanishes on $u(Y)$ and, by
continuity, on $W$. We have $\iota_WB=uA$ and $B\iota_W=0$.
Both $X$ and its closed subspace $W$ have the AP, the latter by HAP.
Thus~\eqref{eq:cyclic}, applied to $B$ and $\iota_W$, gives
\begin{equation}\label{eq:trace-uA}
 \tr_X(uA)=\tr_X(\iota_WB)=\tr_W(B\iota_W)=0.
\end{equation}

On the other hand,
\begin{equation}\label{eq:trace-expansion}
 \begin{aligned}
 \tr_X(uA)
 &=\tr_X(JS)+\tr_X(JVq)+\tr_X(aQS)+\tr_X(aQVq)\\
 &=\tr_X(JS)+\tr_H(QSa)+\tr_H(QVqa)\\
 &=\tr_X(JS)+\tr_H(QR).
 \end{aligned}
\end{equation}
For the vanishing of the second term in the first line, choose a nuclear
representation $V=\sum_k\varphi_k\otimes y_k$ in $\Nuc(X/M,Y)$.
Then
\[
 \tr_X(JVq)=\sum_k\varphi_k(qJy_k)=0.
\]
The series is absolutely convergent, and $qJ=0$. The remaining
two terms are exchanged using~\eqref{eq:cyclic}: the operators
$QS:X\to H$ and $QVq:X\to H$ are nuclear, while $a:H\to X$
is bounded. Since both $X$ and $H$ have the AP,
\[
 \tr_X(aQS)=\tr_H(QSa),\qquad
 \tr_X(aQVq)=\tr_H(QVqa).
\]
Finally, $QSa+QVqa=QR$. Thus
\eqref{eq:trace-expansion}, together with~\eqref{eq:trace-uA}, gives
\eqref{eq:trace-zero}. Notice that no trace on $Y$ or $X/M$ has
been used.
\end{proof}

We apply the proposition to a nuclear witness supplied by
Lemma~\ref{lem:trace-witness}.

\begin{corollary}\label{cor:infinite-codim}
Let $X$ have the HAP and $H$ be a separable Hilbert space. Let
$Y$ be a closed subspace of $X\oplus_2H$ whose horizontal projection
is contained in a closed subspace $M\subset X$ of infinite
codimension. Then $Y$ has the AP.
\end{corollary}

\begin{proof}
Write $J:Y\to X$ and $Q:Y\to H$ for the coordinate maps, and
let $\iota:Y\to E=X\oplus_2H$ be the inclusion. The space $E$
has the AP. If $Y$ failed the AP,
Lemma~\ref{lem:trace-witness} would give $T\in\Nuc(E,Y)$ with
$T\iota=0$ and $\tr_E(\iota T)=1$. Restricting $T$ to the two
summands gives
\[
 T(x,h)=Sx+Rh,\qquad S\in\Nuc(X,Y),\quad R\in\Nuc(H,Y).
\]
Since $T\iota=0$ and $\iota y=(Jy,Qy)$, we obtain
$SJ+RQ=0$. To verify the diagonal trace identity, write
$T=\sum_k e_k^*\otimes y_k$ and decompose
$e_k^*(x,h)=x_k^*(x)+h_k^*(h)$. The induced representations of $S$
and $R$ are nuclear into $Y$, and absolute convergence gives
\[
 \begin{aligned}
 1=\tr_E(\iota T)
 &=\sum_k\bigl(x_k^*(Jy_k)+h_k^*(Qy_k)\bigr)\\
 &=\tr_X(JS)+\tr_H(QR).
 \end{aligned}
\]
This contradicts Proposition~\ref{prop:transfer}.
\end{proof}

\section{A change of coordinates for a closed graph}\label{sec:shear}

To apply Corollary~\ref{cor:infinite-codim} to an arbitrary subspace,
we need a change of coordinates after which the projection onto $X$
has range in a subspace of infinite codimension. The next lemma
gives the required construction for a closed graph. It is convenient
to state it for operators with a reflexive codomain.
The proof uses the annihilator of the graph to obtain a normalized
weak-star null sequence in $X^*$. We select a subsequence and
construct bounded vectors biorthogonal to it. A summable sequence
of functionals on the codomain then defines the nuclear perturbation.
Reflexivity is used to obtain norm density of the domain of the
adjoint in the dual of the codomain.
For an operator $G:D(G)\subset X\to Z$, its graph is
$\Graph G=\{(x,Gx):x\in D(G)\}$. The operator is closed if its graph
is closed in $X\oplus_2Z$, and densely defined if $D(G)$ is dense in $X$.

\begin{lemma}\label{lem:shear}
Let $X$ be a Banach space, let $Z$ be a reflexive Banach space, and let
$G:D(G)\subset X\to Z$ be a closed, densely defined, unbounded linear
operator. For every $\varepsilon>0$, there are $C\in\Nuc(Z,X)$ with
$\nu(C)<\varepsilon$ and a closed subspace $M\subset X$ of infinite
codimension such that
\begin{equation}\label{eq:shear}
 x-CGx\in M\qquad(x\in D(G)).
\end{equation}
\end{lemma}

\begin{proof}
We first find a sequence $z_n^*\in D(G^*)$ such that
$\norm{G^*z_n^*}=1$ and $\norm{z_n^*}\to0$.
Define $D(G^*)$ to consist of the $z^*\in Z^*$ for which $z^*G$ extends
to a bounded functional on $X$, and let $G^*z^*$ be that extension.
It is unique because $D(G)$ is dense. The annihilator of the graph is
\[
 (\Graph G)^\perp
 =\{(-G^*z^*,z^*):z^*\in D(G^*)\}.
\]
Indeed, a pair $(x^*,z^*)$ annihilates the graph precisely when
$z^*G=-x^*$ on $D(G)$.

We claim that $D(G^*)$ is norm dense in $Z^*$. Otherwise,
Hahn--Banach and reflexivity give a nonzero $z\in Z$ annihilating
$D(G^*)$. The displayed identity implies that $(0,z)$ annihilates
$(\Graph G)^\perp$. Since $\Graph G$ is closed, Hahn--Banach
separation gives $(0,z)\in\Graph G$, which forces $z=G0=0$.
This proves the claim.

The map $G^*:D(G^*)\to X^*$ is unbounded when its domain carries the
$Z^*$ norm. Otherwise there would be $K<\infty$ such that
$\norm{G^*z^*}\le K\norm{z^*}$, and the density just proved would give
\[
 \norm{Gx}
 =\sup_{\substack{z^*\in D(G^*)\\\norm{z^*}\le1}}|z^*(Gx)|
 \le K\norm{x}\qquad(x\in D(G)).
\]
Thus $G$ would be bounded on its domain, contrary to the hypothesis.
We may therefore
choose $z_n^*\in D(G^*)$ with
\[
 \norm{G^*z_n^*}=1,\qquad\norm{z_n^*}\longrightarrow0.
\]

Put $f_n=G^*z_n^*$. This is a normalized weak-star null sequence in
$X^*$. In fact, for $x_0\in D(G)$, $f_n(x_0)=z_n^*(Gx_0)\to0$.
For arbitrary $x\in X$,
\[
 |f_n(x)|\le\norm{x-x_0}+|z_n^*(Gx_0)|.
\]
Taking the limit superior and then using density of $D(G)$ shows
that $f_n(x)\to0$ for every $x\in X$.

\smallskip
We next select a subsequence $(\varphi_i)$ of $(f_n)$ and a bounded
sequence $(z_i)$ in $X$ such that $\varphi_i(z_j)=\delta_{ij}$.
For this we use the following elementary observation.
For every finite-dimensional subspace $F\subset X^*$,
\[
 \liminf_n\dist(f_n,F)\ge\frac12.
\]
Otherwise, for some $r<1/2$, there would be infinitely many $n$ and
$e_n\in F$ with $\norm{f_n-e_n}<r$. Since $\norm{e_n}\le1+r$,
we may pass to a subsequence on which $e_n$ converges in norm to
$e\in F$. Now $f_n-e_n\to-e$ weak-star, so weak-star lower
semicontinuity gives $\norm{e}\le r$. But
$1\le\norm{f_n-e_n}+\norm{e_n}$ gives in the limit
$1\le r+\norm{e}\le2r$, a contradiction.

Inductively choose strictly increasing indices $n_i$ and vectors
$v_i\in X$. At stage $i$, write $\varphi_j=f_{n_j}$ for the already
chosen functionals, and choose $n_i$ so large that, with
$\varphi_i=f_{n_i}$ and $\psi_i=z_{n_i}^*$,
\[
 \begin{gathered}
 \dist\bigl(\varphi_i,\spn\{\varphi_j:j<i\}\bigr)>\frac13,
 \qquad |\varphi_i(v_j)|\le2^{-i-j-3}\quad(j<i),\\
 \norm{\psi_i}<\varepsilon2^{-i}/10.
 \end{gathered}
\]
These choices are possible by the distance estimate, weak-star
convergence of $(f_n)$, and $\norm{z_n^*}\to0$.

Let $K_i=\bigcap_{j<i}\ker\varphi_j$, with $K_1=X$.
Its annihilator is
$K_i^\perp=\spn\{\varphi_j:j<i\}$, a finite-dimensional subspace of
$X^*$. The restriction map $X^*\to K_i^*$ is a quotient map by
Hahn--Banach, so
\[
 \norm{\varphi_i|_{K_i}}
 =\dist(\varphi_i,K_i^\perp)>\frac13.
\]
Choose a unit vector $w_i\in K_i$ with $|\varphi_i(w_i)|>1/4$,
which is possible because the restriction norm exceeds $1/3$.
Set $v_i=w_i/\varphi_i(w_i)$. Then
\[
 \varphi_i(v_i)=1,\qquad\norm{v_i}<4,
 \qquad\varphi_j(v_i)=0\quad(j<i).
\]

Write $D=(\varphi_i(v_j))_{i,j\ge1}=I+L$. The matrix $L$ is
strictly lower triangular, and
\[
 \sum_{i>j}|L_{ij}|
 \le\sum_{i>j}2^{-i-j-3}=2^{-2j-3}.
\]
Thus $L$ defines an operator on $\ell_1$ with
\[
 \norm{L}_{\ell_1\to\ell_1}\le\frac1{32}<\frac18.
\]
Consequently $D$ is invertible on $\ell_1$, with
$\norm{D^{-1}}\le8/7$.

Define $V:\ell_1\to X$ by
$Vb=\sum_j b_jv_j$, which is well defined with $\norm{V}\le4$.
For every $b\in\ell_1$ and each fixed $i$, absolute convergence gives
$\varphi_i(Vb)=\sum_jb_j\varphi_i(v_j)=(Db)_i$. Let $(e_j)$ be the standard unit vector basis of $\ell_1$. If
$z_j=VD^{-1}e_j$, then
\begin{equation}\label{eq:biorthogonal}
 \varphi_i(z_j)=\delta_{ij},\qquad
 \sup_j\norm{z_j}\le\frac{32}{7}<5,
 \qquad G^*\psi_i=\varphi_i.
\end{equation}
Indeed, $\varphi_i(z_j)=(DD^{-1}e_j)_i=\delta_{ij}$, and
$\norm{z_j}\le\norm{V}\norm{D^{-1}}\le32/7$.

\smallskip
We can now define the required operator by
\begin{equation}\label{eq:shear-definition}
 Cz=\sum_i\psi_i(z)z_i,\qquad M=\bigcap_i\ker\varphi_i.
\end{equation}
Then
\[
 \sum_i\norm{\psi_i}\norm{z_i}
 \le\frac{32}{7}\sum_i\norm{\psi_i}
 <\frac{16}{35}\varepsilon<\varepsilon.
\]
Hence $C\in\Nuc(Z,X)$ and $\nu(C)<\varepsilon$. The space $M$ is
closed as an intersection of kernels of continuous functionals.
The cosets $z_i+M$ are linearly independent: applying
$\varphi_j$ to any finite linear combination lying in $M$ shows that
its $j$th coefficient is zero. Thus $M$ has infinite codimension.
Finally, if $x\in D(G)$, then $\varphi_n(x)=\psi_n(Gx)$.
The absolutely convergent series defining $C$ and
\eqref{eq:biorthogonal} give
\[
 \varphi_n(x-CGx)
 =\varphi_n(x)-\sum_i\psi_i(Gx)\varphi_n(z_i)
 =\varphi_n(x)-\psi_n(Gx)=0.
\]
Hence $x-CGx\in M$, as required.
\end{proof}

The small nuclear norm in the lemma gives a corresponding control
of the change of coordinates. Indeed, let
\[
 U:X\oplus_2Z\to X\oplus_2Z,\qquad U(x,z)=(x-Cz,z).
\]
Then $U^{-1}(x,z)=(x+Cz,z)$, and
\[
 \norm{U-I}\le\norm{C}<\varepsilon,\qquad
 \norm{U^{-1}-I}\le\norm{C}<\varepsilon.
\]
Thus the graph can be put in the required position by an isomorphism
arbitrarily close to the identity. This conclusion uses neither an
approximation property of $X$ nor a Hilbert structure on $Z$.
The Hilbert structure will enter separately when we split off the
vertical subspace in the proof below.

\section{Proof of the main theorem}\label{sec:completion}

\begin{proof}[Proof of Theorem~\ref{thm:main}]
Suppose that a closed subspace $Y\subset X\oplus_2\ell_2$ fails the AP.
By Lemma~\ref{lem:separable-reduction}, we may assume that $Y$ is
separable.

Write $J:Y\to X$ and $Q:Y\to\ell_2$ for the coordinate maps.
Replacing $X$ by $\cl{J(Y)}$, which still has the HAP, we may also
assume that $X$ is separable and $J(Y)$ is dense in $X$. The space
$Y$ remains closed in the smaller direct sum.
Set
\[
 N=\{h\in\ell_2:(0,h)\in Y\}.
\]
The subspace $N$ is closed. Let $H$ be the orthogonal complement of
$N$ in $\ell_2$. If $P_N$ is the orthogonal projection onto $N$,
the map $y\mapsto(0,P_NQy)$ has its values in $Y$ and is a bounded
projection onto $\{0\}\oplus N$. Its kernel is
\[
 Y_0=\{y\in Y:P_NQy=0\}.
\]
By orthogonality of $N$ and $H$, $Y$ is isometrically isomorphic to
$Y_0\oplus_2N$. Since $N$ has the AP
and the AP passes to finite direct sums, $Y_0$ fails the AP. It is
a closed subspace of $X\oplus_2H$. Write $Q_0:Y_0\to H$ for
the restriction of $Q$ with this codomain.

The restriction $J_0=J|_{Y_0}$ is injective: a vector in its kernel is
vertical and has its second coordinate in $N\cap H=\{0\}$.
Also $J_0(Y_0)=J(Y)$, since subtracting $(0,P_NQy)$ does not change
the horizontal coordinate. Thus $J_0(Y_0)$ is dense in the new $X$.
If $J_0(Y_0)$ were closed, $J_0$ would be a bounded bijection
from $Y_0$ onto $X$. The bounded inverse theorem would then give
$Y_0\simeq X$, contrary to the fact that $Y_0$ fails the AP.
Hence $J_0(Y_0)$ is not closed.

Define
\[
 D(G)=J_0(Y_0),\qquad G(J_0y)=Q_0y\quad(y\in Y_0).
\]
Since $J_0$ is injective with dense range, $G$ is well defined and
densely defined. Its graph is $Y_0$, so $G$ is closed. If $G$ were
bounded on $D(G)$ with respect to the norm of $X$, it would extend
to a bounded operator $\widetilde G:X\to H$. For any $x\in X$,
choose $x_k\in D(G)$ with $x_k\to x$. Then
$(x_k,Gx_k)\to(x,\widetilde Gx)$. Closedness of the graph would
imply $x\in D(G)$, and hence $D(G)=X$, contrary to the fact that
$J_0(Y_0)$ is not closed. Thus $G$ is unbounded.

Apply Lemma~\ref{lem:shear} to obtain a nuclear map $C:H\to X$ and
a closed subspace $M\subset X$ of infinite codimension satisfying
\eqref{eq:shear}. The map
\[
 U:X\oplus_2H\to X\oplus_2H,\qquad U(x,h)=(x-Ch,h)
\]
is a bounded isomorphism with inverse $U^{-1}(x,h)=(x+Ch,h)$.
Hence $Y_1=U(Y_0)$ is closed and fails the AP. Its elements have the form $(x-CGx,Gx)$ for $x\in D(G)$,
so their horizontal coordinates belong to $M$.
Corollary~\ref{cor:infinite-codim} implies that $Y_1$ has the AP, a
contradiction.
\end{proof}

\begin{remark}\label{rem:assumptions}
The proof does not require $X$ to be reflexive, to have a basis, or to
contain an isomorphic copy of $\ell_2$. In Lemma~\ref{lem:shear},
reflexivity of $Z$ ensures norm density of $D(G^*)$ in $Z^*$.
In Proposition~\ref{prop:transfer}, $u$ need not be injective and its
range need not be closed. Nuclearity into $\cl{u(Y)}$ follows by bounded
composition.

The Hilbert structure has two further roles in the main theorem.
Orthogonal complementation separates the vertical subspace $N$ in the
last proof. In Lemma~\ref{lem:factorization}, orthogonal decompositions
give Hilbert blocks with coordinate projections of norm at most one
and nuclear-norm tail estimates for $R$. These blocks can then be matched
with the uniformly complemented Euclidean sections of
Lemma~\ref{lem:sections}. Thus Lemma~\ref{lem:shear} by itself does not extend the theorem to
an arbitrary reflexive HAP summand.

The conclusion does not settle Question~9 of~\cite{JS}, concerning
$\ell_2(X)$. Iterating the theorem adds Hilbert summands, but does
not give the HAP for $X\oplus X$. Thus such iteration does not
address even the finite sums of copies of $X$ that would arise in
an approach to $\ell_2(X)$ by finite coordinate subspaces.

For comparison, Anisca and Chlebovec~\cite[Theorem~1.1]{AC} proved
that if $E=\ell_2(Z)$ has cotype~2 for a Banach space $Z$, and
there are $\alpha>0$, $\beta>1$, and $n_0\in\mathbb N$ such that
\[
 d_n(E)\ge\alpha(\log n)^\beta\qquad(n\ge n_0),
\]
then $E$ has a subspace without the AP. Their result concerns a
class of countable Hilbertian sums under a quantitative hypothesis
on local Euclidean distances. The present theorem concerns the
addition of a Hilbert space to a single HAP space and requires no
such quantitative hypothesis.
\end{remark}

The extension to a nonseparable Hilbert summand follows from the
separable reduction in Lemma~\ref{lem:separable-reduction}.

\begin{corollary}\label{cor:arbitrary-hilbert}
If $X$ has the HAP and $H$ is any Hilbert space, then $X\oplus_2H$ has the HAP.
\end{corollary}

\begin{proof}
Suppose that a closed subspace of $X\oplus_2H$ fails the AP.
By Lemma~\ref{lem:separable-reduction}, it contains a separable
closed subspace $Y_s$ without the AP. Let $Q:Y_s\to H$ be the
coordinate map and put $H_0=\cl{Q(Y_s)}$. Then $H_0$ is a
separable Hilbert space and $Y_s\subset X\oplus_2H_0$. An
isometric embedding $j:H_0\to\ell_2$ maps $Y_s$, by
$I_X\oplus j$, onto a closed subspace of $X\oplus_2\ell_2$
without the AP. This contradicts Theorem~\ref{thm:main}.
\end{proof}

In the reflexive case, every closed subspace of the direct sum has
the MAP.

\begin{corollary}\label{cor:hmap}
Let $X$ be a reflexive Banach space with the HAP, and let $H$ be
any Hilbert space. Then $X\oplus_2H$ has the HMAP.
\end{corollary}

\begin{proof}
By Corollary~\ref{cor:arbitrary-hilbert}, $X\oplus_2H$ has the HAP.
It is reflexive, since both $X$ and $H$ are reflexive. Let $Y$ be
any closed subspace of $X\oplus_2H$. Then $Y$ is reflexive and has
the AP. By Grothendieck's theorem
(see~\cite[Theorem~1.e.15 and pp.~40--41]{LT}), $Y$ has the MAP.
Since $Y$ was arbitrary, $X\oplus_2H$ has the HMAP and, in
particular, the HBAP.
\end{proof}

\section*{AI Declaration}
During the preparation of this work, the authors used GPT-6 Astra to
assist with language editing and the refinement of mathematical
arguments. The authors take full responsibility for the content of
this paper.

\section*{Competing interests}
On behalf of all authors, the corresponding author states that
there is no conflict of interest.

\end{document}